\documentclass{amsart}
\usepackage{amsthm, amscd, amssymb, amsfonts,dsfont}
\usepackage{mathrsfs}
\usepackage[latin1]{inputenc}
\usepackage{graphicx}
\usepackage{enumerate}
\usepackage{esint}
\usepackage{srcltx}
\usepackage{hyperref}

\newtheorem{thm}{Theorem}

\newtheorem{lem}{Lemma}
\newtheorem{prop}{Proposition}
\theoremstyle{remark}
\newtheorem{rem}{Remark}

\newcommand{\N}{\mathds N}
\newcommand{\R}{\mathds R}
\newcommand{\te}{\theta}
\newcommand{\p}{p(\cdot)}
\newcommand{\pri}{p'(\cdot)}
\newcommand{\be}{\beta}
\newcommand{\loc}{\text{\upshape loc}}
\newcommand{\Bro}{\mathcal{B}_\rho}

\begin{document}

\title[Weighted norm inequalities for the maximal functions associated ...]{Weighted norm inequalities for the maximal functions associated to a critical radius function on variable Lebesgue spaces}

\author{A. Cabral}



\thanks{This research is partially supported by grants from Consejo Nacional de Investigaciones Cient\'ificas y T\'ecnicas (CONICET) and Universidad Nacional del Nordeste (UNNE), Argentina.}

\address{Departamento de Matem\'atica -- Facultad de Ciencias Exactas y Naturales  -- UNNE and Instituto de Modelado e Innovaci\'on Tecnol\'ogica, CONICET--UNNE, Corrientes, Argentina.}

\email{enrique.cabral@comunidad.unne.edu.ar}

\date{}

\begin{abstract}
In this work we obtain boundedness on weighted variable Lebesgue spaces of some maximal functions that come from the localized analysis consi-dering a critical radius function. This analysis appears naturally in the context of the Schr\"odinger operator $\mathcal{L}=-\Delta+V$ in $\R^d$,  where $V$ a non-negative potential which satisfying a some specific reverse H\"older condition. We consider new classes of weights that locally behave as the Muckenhoupt class for Lebesgue spaces with variable exponents considered in~\cite{CU-Fiorenza-Neuge} and actually including them.
\end{abstract}

\maketitle

\section{Introduction and Preliminaries}

A critical radius function $\rho$ is an application that assigns to each $x$ in $\R^d$ a positive number in a way that its variation at different points is somehow controlled by a power of the distance between them. 

More specifically, we call a \emph{critical radius function} to any positive  function $\rho$ with the property that there exist constants $c_\rho$, $N_0\geq1$ such that
\begin{equation}\label{est rho}
 c_\rho^{-1}\rho(x)\bigg(1+\frac{|x-y|}{\rho(x)}\bigg)^{-N_0}\leq\rho(y)\leq c_\rho\,\rho(x)\bigg(1+\frac{|x-y|}{\rho(x)}\bigg)^{\frac{N_0}{N_0+1}},
\end{equation}
for every $x,y\in\R^d$. Clearly, if $\rho$ is such a function, so it is $\beta\rho$ for any $\beta>0$.

The constant functions are a trivial example of critical radius functions. So is function $\frac{1}{1+|x|}$ (see Section 6).

This kind of function appears naturally in the harmonic analysis  related to a Schr{\"o}dinger operator $\mathcal{L}=-\Delta+V$ in $\R^d$, $d\ge 3$, where $V$ is a non-negative, non-identically zero and  satisfies a reverse H\"older condition for $q>d/2$,
\begin{equation*}\label{reverse}
\left(\frac{1}{|B|}\int_B V(y)^q\,dy\right)^{1/q} \leq \frac{C}{|B|}
\int_B V(y)\,dy,
\end{equation*}
for every ball $B\subset \R^d$.

Given a critical radius function $\rho$,  a ball $B(x,r)\subset\R^d$ will be called \emph{critical} if $r=\rho(x)$. We denote by $\Bro$ the family of \emph{sub-critical} balls, i.e.
\begin{equation}\label{def:Bro}
\mathcal{B}_\rho=\{B(x,r): x\in\R^d\ and\ r\le\rho(x)\}.
\end{equation}

One of the operators we are interested in is the local maximal operator $M^\loc_\rho$ defined for $f\in L^1_{\text{loc}}(\R^d)$ as
\begin{equation}\label{max_local}
M^\loc_\rho f(x) = \sup_{r\le\rho(x)}\frac{1}{|B(x,r)|} \int_{B(x,r)}|f(y)|\,dy.
\end{equation}

When $\rho\equiv\infty$, the operator $M^\loc_\rho$ coincides with the maximal Hardy--Littlewood function and in this case the problem for weighted variable exponent was completely
solved in~\cite{CU-Diening-Hasto} and~\cite{CU-Fiorenza-Neuge}. On the other hand,  in the case of a  general $\rho$ and $p$ constant, in \cite{BHS-Classes} the authors proved that $M^\loc_\rho$ is bounded on $L^p(w)$ for $1<p<\infty$ where $w$ belongs to a suitable class larger than the classical $A_p$ Muckenhoupt weights. 

In addition to the localized maximal operator presented in~\eqref{max_local} , we are interested in the following maximal operators
\begin{equation}\label{Max_theta}
M^\te_{\rho} f(x) = \sup_{r>0}\bigg(1+\frac r{\rho(x)}\bigg)^{-\te} \frac{1}{|B(x,r)|} \int_{B(x,r)}|f(y)|\,dy,
\end{equation}
for $\theta>0$ and $f\in L^1_{\text{loc}}(\R^d)$.

The boundedness in $L^p(w)$ with appropriate weights of these operators was duly studied in~\cite{BCH-Extrapol}.

The aim of this paper is to determine the weights for which the operators $M^\loc_\rho$ and $M^\te_{\rho}$, $\theta>0$, preserve  the $L^{\p}(w)$ spaces, known as weighted variable Lebesgue spaces.

The remainder of this paper is organized as follows. In the rest of this section we demonstrate some facts about the critical radius function. In Section 2 we gather together some basic results about variable Lebesgue spaces that will be used later. Section 3 is devoted to present the classes of weights involved in this work and some of their properties that will be useful. In Section 4 we state and 	prove the main results of this work by finding the behavior of the maximal operators we have already talked about, and finally we present some applications to the context of the Schr\"odinger operator in Section 5.

Throughout this paper, unless otherwise indicated, we will use $C$ and $c$ to denote constants, which are not necessarily the same at each occurrence.

\

A useful tool is that whenever $\rho$ satisfies~\eqref{est rho}, it is possible to cover $\R^d$ by a countable family of critical balls  in a way that any fixed dilation of the covering has the bounded overlapping property. More precisely,

\begin{prop}[See \cite{DZ-H1}]\label{propo: cubr}
There exists a sequence of points $\{x_k\}_{k\ge1}$ in $\R^d$ such that family  $B_k=B(x_k,\rho(x_k))$ satisfies
\begin{enumerate}[i)]
\item $\displaystyle \cup_k B_k = \R^d$;
\item \label{overlap} for every $\be\ge 1$ there exist constants $C$ and $N_1$ such that, $\sum_k \chi_{\be B_k} \leq C\be^{N_1}$.
\end{enumerate}
\end{prop}

Given a ball $B_0\notin\mathcal{B}_\rho$ we shall also need a particular covering by sub-critical balls with centers inside $B_0$ as the following lemma shows.

\begin{prop}\label{prop:cubB}
Let $B_0=B(x_0,r)$ with $x_0\in\R^d$ and $\rho(x_0)<\,r\le\be\rho(x_0)$, where $\be>1$ . Then, there exists a covering of $B_0$ by balls $\{P_j\}_{j=1}^N$ such that 
\begin{enumerate}[i)]
\item\label{prop} $\frac12P_j\cap\frac12P_k=\emptyset$, for $k\neq j$;
\item $\sum_{i=1}^N \chi_{P_i} \leq C$, where $C$ depends only on $d$ and $N$ depends on $\be$  and the constants in \eqref{est rho}.
\end{enumerate}
\end{prop}
\begin{proof}
First observe that for $x\in2B_0$, by~\eqref{est rho} we have
\begin{equation}\label{delta_cero}
\rho(x) \ge c_{\rho}^{-1}\left(1+\frac{2r}{\rho(x_0)}\right)^{-N_0}\rho(x_0) = \delta_0> 0.
\end{equation}

Consider the family of sets
\begin{equation*}
\mathcal{F}=\left\{S\subset B_0:\ B\Big(x,\frac{\delta_0}{8}\Big)\cap B\Big(y,\frac{\delta_0}{8}\Big)=\emptyset,\ \ \forall\, x,y \in S,\ x\neq y\right\}.
\end{equation*}

It is clear that $\mathcal{F}\neq\emptyset$ since $\{x_0\}\in\mathcal{F}$. Observe that if $\mathcal C$ be a chain in $\mathcal{F}$ endowed with the order of inclusion, then $V=\cup_{S\in\mathcal C} S$ is an upper bound of $\mathcal C$. Therefore, there exists a maximal element $S_{\max}$ in $\mathcal{F}$. 

Denote by $x_1,x_2,\ldots,x_N$ the elements of $S_{\max}$ and $P_i=B(x_i,\frac{\delta_0}{4})$. It is clear then that  $\frac12P_j\cap\frac12P_k=\emptyset$, for $k\neq j$. Also, property~\eqref{prop}  implies that $N$ is bounded by $C_1\beta^{d(N_0+1)}$, with $C_1$ depending only on the dimension $d$ and $c_\rho$.

We shall see that $B_0\subset \cup_{i=1}^N P_i$ and the overlapping of the balls $P_i$, $i=1,\ldots,N$, is bounded.
	
Suppose there exists $y\in B_0$ such that $y\notin \cup_{i=1}^N P_i$, which means $|y-x_i|\ge\frac{ \delta_0}{4}$, $i=1,\ldots,N$. Now  see that $B(y, \delta_0/8)\cap B(x_i, \delta_0/8)$ is empty. In fact, suppose $z\in B(y, \delta_0/8)\cap B(x_i, \delta_0/8)$, then 
\begin{equation*}
|y-x_i|\le |y-z| + |z-x_i| < \frac{ \delta_0}{8}+\frac{ \delta_0}{8}=\frac{ \delta_0}{4},
\end{equation*} 
which is a contradiction. So, $S_{\max}\cup \{y\}$ belongs to $\mathcal{F}$  and this means the contradiction that $S_{\max}$ is not a maximal element of $\mathcal{F}$.
	
Now we see that the overlapping of the balls $P_i$, $i=1,\ldots,N$, is finite and depends only on $d$.
	
Suppose that $m$ is such that $\cap_{i=1}^m B(y_i, \delta_0/4)\neq\emptyset$ for some points $y_i\in S_{\max}$. Then, we have
\begin{equation*}
\bigcup_{i=1}^m B\Big(y_i,\frac{ \delta_0}{8}\Big)\subset B\Big(y_1,\frac{5}{8} \delta_0\Big).
\end{equation*}
Now, due to the fact that the balls $\{B(y_i, \delta_0/8)\}_{i=1}^m$ are disjoints we conclude
\begin{equation*}
\begin{split}
m |B(y_1, \delta_0/8)| &= \sum_{i=1}^{m} |B(y_i, \delta_0/8)| \\
& =|\cup_{i=1}^m B(y_i, \delta_0/8)|\\
& \le |5B(y_1, \delta_0/8)|,
\end{split}
\end{equation*}
thus $m\leq5^d$.
\end{proof}

\section{Variable Lebesgue Spaces}
We begin with some definitions related to the variable Lebesgue spaces. Let $\p:\R^d\rightarrow[1,\infty]$ be a measurable function. Given a measurable set $A\subset\R^d$ we define
\begin{equation*}
p^{-}_{A}:=\text{ess}\inf_{x\in A}p(x),\hspace{1cm} p^{+}_{A}:=\text{ess}\sup_{x\in A}p(x).
\end{equation*}

For simplicity we denote $p^{+}=p^{+}_{\R^d}$ and $p^{-}=p^{-}_{\R^d}$. Given $\p$, the conjugate exponent $\pri$ is defined pointwise
\begin{equation*}
\frac{1}{p(x)}+\frac{1}{p'(x)}=1,
\end{equation*}
where we let $p'(x)=\infty$ if $p(x)=1$.

By $\mathcal{P}(\R^d)$ we will designate  the collection of all measurable functions $\p:\R^d\rightarrow[1,\infty]$ and by $\mathcal{P}^*(\R^d)$ the set of $p\in\mathcal{P}(\R^d)$ such that $ p^{+}<\infty$.

Given  $p\in\mathcal{P}^*(\R^d)$, we say that a measurable function $f$ belongs to $L^{\p}(\R^d)$ if for some $\lambda>0$
\begin{equation*}
\varrho(f/\lambda)=\int_{\R^d}\left(\frac{|f(x)|}{\lambda}\right)^{p(x)}dx<\infty.
\end{equation*}
A Luxemburg norm can be defined in $L^{\p}(\R^d)$ by taking
\begin{equation*}
\|f\|_{L^{\p}(\R^d)}=\inf\{\lambda>0:\varrho(f/\lambda)\le1\}.
\end{equation*}

These spaces are special cases of Musieliak-Orlicz spaces (see \cite{Musielak}), and generalize the classical Lebesgue spaces. For more information see,  for example \cite{KR, Libro:CU-Fiorenza, Libro:D-H}.

We will indicate with $L^{\p}_{\loc}(\R^d)$ the space of functions $f$ such that $f\chi_{B}\in L^{\p}(\R^d)$ for every ball $B\subset\R^d$.

The following conditions on the exponent arise in connection with the boundedness of the Hardy--Littlewood maximal operator $M$ in $L^{\p}(\R^d)$ (see, for example, \cite{Diening}, \cite{Libro:CU-Fiorenza} or~\cite{Libro:D-H}). We will say that $p\in\mathcal{P}^{\log}(\R^d)$ if $p\in\mathcal{P}^*(\R^d)$ and if there are constants $C>0$ and $p_\infty$ such that
\begin{equation}
|p(x)-p(y)|\le \frac{C}{\log(e+1/|x-y|)},\hspace{1cm}\text{for every}~x,y\in\R^d,
\end{equation}
and
\begin{equation}
|p(x)-p_\infty|\le \frac{C}{\log(e+|x|)},\hspace{1cm}\text{with}~x\in\R^d.
\end{equation}
The functions in $\mathcal{P}^{\log}(\R^d)$ are generally called log--H\"older continuous.

With this definition of $\|\cdot\|_{\p}$ we have the following generalization of H\"older's inequality and an equivalent expression for the norm (see, for example,~\cite{Libro:D-H}).

\begin{lem}\label{lem: Holder}
Let $p\in\mathcal{P}(\R^d)$. Then, for every $f\in L^{\p}(\R^d)$ and $g\in L^{p'(\cdot)}(\R^d)$, $fg\in L^1(\R^d)$ and
\begin{equation*}
\int_{\R^d}|f(x)||g(x)|\,dx\le2\|f\|_{\p}\|g\|_{\pri}.
\end{equation*}
\end{lem}

\begin{lem}\label{lem: formula norma}
Let $p\in\mathcal{P}(\R^d)$. Then,
\begin{equation*}
\frac12\|f\|_{\p}\le\sup_{\|g\|_{\pri}\le1}\ \int_{\R^d}|f||g|\,dx\le2\|f\|_{\p}.
\end{equation*}
\end{lem}

For a weight $w$ means a locally integrable function such that $0<w(x)<\infty$ almost everywhere. 

Given a weight $w$ and $p\in\mathcal{P}(\R^d)$, there are two possible ways to define the weighted variable Lebesgue space. One way is to treat $w\,dx$ as a measure, and define the weighted variable Lebesgue space $L^{\p}(w)$ with regard to this measure. This approach is in one respect a natural generalization of $L^p(w)$ and was adopted by Diening and H\"ast\"o in~\cite{Diening-Hasto}.

We adopt a different approach, following the spirit of~\cite{CU-Fiorenza-Neuge}, we replace the weight $w$ by $w^p$, that is, we considering the weight $w$ acts as a multiplier. Given a weight $w$ and $p\in\mathcal{P}(\R^d)$, we define the weighted variable Lebesgue space $L^{\p}(w)$ to be the set of all measurable functions $f$ such that $fw\in L^{\p}$, and we write $\|f\|_{L^{\p}(w)}=\|f\|_{\p,w}=\|fw\|_{\p}$. Thus, we say that an operator $T$ is bounded on
$L^{\p}(w)$ if $\|(Tf)w\|_{\p}\le C\|fw\|_{\p}$ for all $f\in L^{\p}(w)$.

In the rest of this section we display a simple and convenient method to pass from local to global results in variable exponent spaces. The idea is simply to generalize the following two properties of the Lebesgue norm,
\begin{equation*}
\|f\|_{L^q}=\left(\sum_{k}\|\chi_{\Omega_k} f\|^q_{L^q}\right)^{1/q},
\end{equation*}
and
\begin{equation*}
\sum_{k}\|\chi_{\Omega_k} f\|_{L^q}\|\chi_{\Omega_k} g\|_{L^{q'}}\le\|f\|_{L^q}\|g\|_{L^{q'}},
\end{equation*}
for a partition of $\R^d$ into measurable sets $(\Omega_k)_k$ and $1\le q<\infty$.

This technique was studied and developed in~\cite{Libro:D-H} (see also~\cite{H}).
Although these estimates cannot be generalized for arbitrary choices of the $(\Omega_k)_k$ sets, as shown in Theorem 7.3.22 they hold as long as $p\in\mathcal{P}^{\log}(\R^d)$  and when $(\Omega_k)_k$ is a locally $N$--finite family of balls (or cubes).

A family \ $\mathcal{U}$ of measurable sets $U\subset\R^d$ is called to be locally $N$--finite, if 
\begin{equation*}
\sum_{U\in\mathcal{U}}\chi_U\le N
\end{equation*}
almost everywhere in $\R^d$.

\begin{rem}\label{n-familias}
It is clear from~Proposition~\ref{propo: cubr} that, given $\beta\ge1$, the family of balls $\mathcal{B}=\{\be B_k\}_{k\in\N}$ is locally $N$--finite with $N=C\be^{N_1}$.
\end{rem}
 
The following theorem is nothing but a deep insight of the proof of Theorem~7.3.22 en~\cite{Libro:D-H}. Here for the sake of convenience, we supply the proof.

\begin{thm}\label{teo:H}
Let $p\in\mathcal{P}^{\log}(\R^d)$ and $\mathcal{B}$ a  locally $N$--finite family of balls. Then,
\begin{equation*}
\sum_{B\in\mathcal{B}}\|\chi_B f\|_{\p}\|\chi_B g\|_{\pri}\le CN^2\|f\|_{\p}\|g\|_{\pri},
\end{equation*}
for every $f\in L^{\p}_{\loc}(\R^d)$ and $g\in L^{\pri}_{\loc}(\R^d)$.
\end{thm}
\begin{proof}
Since $p\in\mathcal{P}^{\log}(\R^d)$, it is easy to see that $p'\in\mathcal{P}^{\log}(\R^d)$. Further, $p\in\mathcal{P}^{\log}(\R^d)$ also implies that
\begin{equation}\label{equiv}
|B|\approx\|\chi_B\|_{\p}\|\chi_B\|_{\pri},
\end{equation}
for every ball $B\subset\R^d$. 
 
Consider 	        	 	
\begin{equation*}
T_{\p}^\mathcal{B}f=\sum_{B\in\mathcal{B}}\chi_{B}\frac{\|f\chi_B\|_{\p}}{\|\chi_B\|_{\p}}.
\end{equation*}

Hence, from Corollary~7.3.21 in~\cite{Libro:D-H} (and its proof), it follows that
\begin{equation*}
\|T_{\p}^\mathcal{B}\,f\|_{\p}\le CN\|f\|_{\p}
\hspace{0.3cm}\text{and}\hspace{0.3cm}
\|T_{\pri}^\mathcal{B}\,g\|_{\pri}\le CN\|g\|_{\pri},
\end{equation*}	
for every $f\in L^{\p}_{\loc}(\R^d)$ and $g\in L^{\pri}_{\loc}(\R^d)$. 

Therefore, from \eqref{equiv} and Lemma~\ref{lem: Holder} we get
\begin{equation*}
\begin{split}
\sum_{B\in\mathcal{B}}\|\chi_B f\|_{\p}\|\chi_B\,g\|_{\pri}&\le C\sum_{B\in\mathcal{B}}|B|\frac{\|\chi_B f\|_{\p}}{\|\chi_B\|_{\p}}\frac{\|\chi_B\,g\|_{\pri}}{\|\chi_B\|_{\pri}}\\
&\le C\int_{\R^d}\sum_{B\in\mathcal{B}}\chi_{B}\frac{\|\chi_B f\|_{\p}}{\|\chi_B\|_{\p}}\frac{\|\chi_B\,g\|_{\pri}}{\|\chi_B\|_{\pri}}\,dx\\
&\le C\int_{\R^d}T_{\p}^\mathcal{B}\,f(x)\ T_{\pri}^\mathcal{B}\,g(x)\ dx\\
&\le C\|T_{\p}^\mathcal{B}\,f\|_{\p}\,\|T_{\pri}^\mathcal{B}\,g\|_{\pri}\\
&\le CN^2\|f\|_{\p}\|g\|_{\pri}.
\end{split}
\end{equation*}
\end{proof}

Let us now give a result that will be very useful in what follows.

\begin{thm}\label{teo: local-global}
Let $p\in\mathcal{P}^{\log}(\R^d)$ and a weight $w$. Then, for every locally $N$--finite family of balls $\mathcal{B}$ we have
\begin{equation*}
\bigg\|\sum_{B\in\mathcal{B}}\chi_{B}\frac{\|f\chi_B\|_{\p,w}}{\|\chi_B\|_{\p,w}}\bigg\|_{\p,w}	
\le CN^2\|f\|_{\p,w},
\end{equation*}
for every $f\in L^{\p}(w)$. Moreover,
\begin{equation*}
\bigg\|\sum_{B\in\mathcal{B}}\chi_{B}f \bigg\|_{\p,w}\approx
\bigg\|\sum_{B\in\mathcal{B}}\chi_{B}\frac{\|f\chi_B\|_{\p,w}}{\|\chi_B\|_{\p,w}}\bigg\|_{\p,w}.
\end{equation*}

\end{thm}
\begin{proof}
For convenience, we shall denote
\begin{equation*}
T_{\p,w}^{\mathcal{B}}f=\sum_{B\in\mathcal{B}}\chi_{B}\frac{\|f\chi_B\|_{\p,w}}{\|\chi_B\|_{\p,w}}.
\end{equation*}
	
Using the norm conjugate formula, H\"older's inequality  and Theorem~\ref{teo:H} we estimate
\begin{equation*}
\begin{split}
\bigg\|\sum_{B\in\mathcal{B}}\chi_{B}f\bigg\|_{\p,w}
&\le2\ \sup_{\|g\|_{\pri}\le1}\ \int_{\R^d}\bigg(\sum_{B\in\mathcal{B}}\chi_{B}|f|\bigg)w\,|g|\,dx\\
&\le2\ \sup_{\|g\|_{\pri}\le1}\ \sum_{B\in\mathcal{B}}\int_{B} w|f||g|\,dx\\
&\le4\ \sup_{\|g\|_{\pri}\le1}\ \sum_{B\in\mathcal{B}}
\|\chi_B wf\|_{\p} \|\chi_B g\|_{\pri}\\
&=4\ \sup_{\|g\|_{\pri}\le1}\ \sum_{B\in\mathcal{B}}
\bigg\|\chi_B w\frac{\|\chi_B f\|_{\p,w}}{\|\chi_B\|_{\p,w}}\bigg\|_{\p} \|\chi_B g\|_{\pri}\\
&\le4\ \sup_{\|g\|_{\pri}\le1}\ \sum_{B\in\mathcal{B}}
\bigg\|\chi_B\,w\, T_{\p,w}^{\mathcal{B}}f\bigg\|_{\p} \|\chi_B g\|_{\pri}\\
&\le4CN^2\ \sup_{\|g\|_{\pri}\le1}\|w\,T_{\p,w}^{\mathcal{B}}f\|_{\p}\|g\|_{\pri}\\
&=4CN^2\ \|w\,T_{\p,w}^{\mathcal{B}}f\|_{\p}\sup_{\|g\|_{\pri}\le1}\|g\|_{\pri}\\
&\le4CN^2\bigg\|\sum_{B\in\mathcal{B}}\chi_{B}\frac{\|f\chi_B\|_{\p,w}}{\|\chi_B\|_{\p,w}}\bigg\|_{\p,w}.
\end{split}
\end{equation*}	
	
Using the same tools, we also determine
\begin{equation*}
\begin{split}
\bigg\|\sum_{B\in\mathcal{B}}\chi_{B}\frac{\|f\chi_B\|_{\p,w}}{\|\chi_B\|_{\p,w}}\bigg\|_{\p,w}
&\le2\ \sup_{\|g\|_{\pri}\le1}\ \int_{\R^d}\bigg(\sum_{B\in\mathcal{B}}\chi_{B}
\frac{\|f\chi_B\|_{\p,w}}{\|\chi_B\|_{\p,w}}\bigg)w|g|\,dx\\
&\le2\ \sup_{\|g\|_{\pri}\le1}\ \sum_{B\in\mathcal{B}}
\frac{\|f\chi_B\|_{\p,w}}{\|\chi_B\|_{\p,w}}\,\int_{B}|g|w\,dx\\
&\le4\ \sup_{\|g\|_{\pri}\le1}\ \sum_{B\in\mathcal{B}}
\frac{\|f\chi_B\|_{\p,w}}{\|\chi_B\|_{\p,w}}\|\chi_Bg\|_{\pri}\|\chi_B w\|_{\p}\\
&=4\ \sup_{\|g\|_{\pri}\le1}\ \sum_{B\in\mathcal{B}}
\|wf\chi_B\|_{\p}\|\chi_Bg\|_{\pri}\\
&\le4CN^2\sup_{\|g\|_{\pri}\le1}\|wf\|_{\p}\ \|g\|_{\pri}\\
&=4CN^2\|f\|_{\p,w}\ \sup_{\|g\|_{\pri}\le1}\|g\|_{\pri}\\
&\le4CN^2\|f\|_{\p,w}.
\end{split}
\end{equation*}
\end{proof}

\section{Weights}
To reach  our results we will need to assume that the maximal operator is bounded on weighted variable Lebesgue spaces. The corresponding results regarding the boundedness of the maximal operator on  weighted  variable Lebesgue spaces  have been established in different settings by several authors, see for instance \cite{CU-Diening-Hasto, CU-Fiorenza-Neuge, Libro:D-H}. Following~\cite{CU-Fiorenza-Neuge},  if we replace the weight $w$ by $w^p$ in the definition of weights of Muckenhoupt, the condition $A_p$ takes the form
\begin{equation*}
\sup_{B}|B|^{-1}\|w\chi_B\|_{p}\|w^{-1}\chi_B\|_{p'}<\infty,
\end{equation*}
where the supreme is taken on all balls $B\subset\R^d$ and $1<p<\infty$. This condition results equivalent to the norm inequality
\begin{equation*}
\|Mf\,w\|_p\le C\|f\,w\|_p.
\end{equation*}

This way of defining the weights $A_p$ immediately generalizes to the variable Lebesgue spaces. Given an exponent $p\in\mathcal{P}(\R^d)$ and a weight $w$, we say that $w\in A_{\p}$ if
\begin{equation}\label{Ap punto}
[w]_{A_{\p}}=\sup_{B}|B|^{-1}\|w\chi_B\|_{\p}\|w^{-1}\chi_B\|_{p'(\cdot)}<\infty,
\end{equation}
where the supreme is taken on all balls $B\subset\R^d$.

The following result can be found in~\cite{CU-Fiorenza-Neuge} and \cite{CU-Diening-Hasto}.

\begin{thm}\label{thm: Acot M-w}
Let $p\in\mathcal{P}^{\log}(\R^d)$ with $p^{-}>1$. Then the Hardy--Littlewood maximal operator $M$ is bounded on $L^{\p}(w)$ if and only if $w\in A_{\p}$.
\end{thm}

\begin{rem}\label{rem:loc}
Using the techniques in~\cite{CU-Fiorenza-Neuge} or~\cite{CU-Diening-Hasto} it can be seen that the above theorem is also valid if one considers the maximal operator for functions defined on a (fixed) ball $Q$ and the corresponding weights $A_{\p}(Q)$, where only the balls $B$ such that $|B\cap Q|>0$ must be considered.
\end{rem}

Following the above and the ideas in~\cite{BCH-Extrapol} and~\cite{BHS-Classes} we define some classes of weights associated to a critical radius function $\rho$. Given a critical radius function $\rho$ and $p\in\mathcal{P}(\R^d)$ we shall consider two families of weights. We introduce the $A_{\p}^{\rho,\loc}$ class of weights as those $w$ satisfying~\eqref{Ap punto} for every ball $B\in\Bro$. That is to say, the weights $w$ such that
\begin{equation}\label{Ap_local}
[w]_{A_{\p}^{\rho,\loc}}=\sup_{B\in\Bro}|B|^{-1}\|w\chi_B\|_{\p}\|w^{-1}\chi_B\|_{p'(\cdot)}<\infty.
\end{equation}

Also, given $\te>0$, we will say that $w\in A_{\p}^{\rho,\te}$ if the inequality
\begin{equation}\label{Ap_rho}
\|w\chi_B\|_{\p}\|w^{-1}\chi_B\|_{p'(\cdot)}\le C|B|\left(1+\frac{r}{\rho(x)}\right)^\te,
\end{equation}
holds for all balls $B=B(x,r)$. We denote $A_{\p}^{\rho}=\cup_{\te\ge0}A_{\p}^{\rho,\te}$. It is clear then that $A_{\p}\subset A_{\p}^\rho\subset A_{\p}^{\rho,\loc}$.

In the rest  of this section we evince some basic properties of the weights in $A_{\p}^{\rho,\loc}$ class.

\begin{lem}\label{lem: Ainf}
Let $p\in\mathcal{P}(\R^d)$. If $w\in A_{\p}^{\rho,\loc}$, then there exists  a constant $C$ depending on $[w]_{A_{\p}^{\rho,\loc}}$ such that, given any ball $B\in\mathcal{B}_\rho$ and any measurable set $E\subset B$,
\begin{equation*}
\|\chi_B w\|_{\p}\le C\|\chi_E w\|_{\p}\frac{|B|}{|E|}.
\end{equation*}
\end{lem}

\begin{proof}
Let $B\in\mathcal{B}_\rho$ and $E\subset B$. Then by Hölder's inequality and the $A_{\p}^{\rho,\loc}$ condition, we have 
\begin{equation*}
\begin{split}
|E|&=\int_{\R^d}\chi_{E}(x)w(x)\chi_B(x)w^{-1}(x)\,dx\\
&\le 2\|\chi_{E}w\|_{\p}\|\chi_Bw^{-1}\|_{\pri}\\
&\le 2[w]_{A_{\p}^{\rho,\loc}}|B|\|\chi_{E}w\|_{\p}\|\chi_Bw\|_{\p}^{-1}.
\end{split}
\end{equation*}
\end{proof}

\begin{prop}\label{prop: Aploc beta}
Let  $p\in\mathcal{P}(\R^d)$ and $\be\ge1$, then $A_{\p}^{\beta\rho,\loc}=A_{\p}^{\rho,\loc}$. Moreover, the constant $[w]_{A_{\p}^{\beta\rho,\loc}}$ depends only on the constants in \eqref{est rho}, $\be$ and  $[w]_{A_{\p}^{\rho,\loc}}$.
\end{prop}
\begin{proof}
It is sufficient to check that $A_{\p}^{\rho,\loc}\subset A_{\p}^{\beta\rho,\loc}$ since the other contention always hold for $\beta\ge1$. Let $w\in A_{\p}^{\rho,\loc}$ and $B_0=B(x_0,r)$ with $r\le\beta\rho(x_0)$. It is sufficient to suppose that $r>\rho(x_0)$ 

Consider the ball $B_1=B(x_1,\frac{r}2)$ with $B_1\subset B_0$. We must prove that there exists a constant $\widetilde{C}_1$ such that
\begin{equation*}
\|w\chi_{B_0}\|_{\p}\le \widetilde{C}_1\|w\chi_{B_1}\|_{\p}.
\end{equation*}

Let $\{P_i\}_{i=1}^N$ the covering of $B_0$ provided by  Proposition~\ref{prop:cubB}, with $P_j=B(y_j,\frac{ \delta_0}4)$ and $ \delta_0=c_{\rho}^{-1}\left(1+2r/\rho(x_0)\right)^{-N_0}\rho(x_0)$.

Therefore, for any $j$ and $k$ such that $P_j\cap P_k\neq\emptyset$, we have that $P_j\subset4P_k$. Then, 
\begin{equation*}
\|w\chi_{P_j}\|_{\p}\le\|w\chi_{4P_k}\|_{\p}\le2^{2d+1}[w]_{A_{\p}^{\rho,\loc}}\|w\chi_{P_k}\|_{\p}=C_{w}\|w\chi_{P_k}\|_{\p},
\end{equation*}
where in the last inequality we have used the Lemma~\ref{lem: Ainf}  since the radius of $4P_k$ is $ \delta_0\le\rho(y_k)$, which in turn holds because $y_k\in2B_0$ and \eqref{delta_cero}.

If  $P_{j_1}$ denotes a member of the family such that $x_1\in P_{j_1}$,  we have 

\begin{equation}\label{dup_w}
\|w\chi_{B_0}\|_{\p}\le\sum_{i=1}^N\|w\chi_{P_j}\|_{\p}\le C_{w}^N\sum_{i=1}^N\|w\chi_{P_{j_1}}\|_{\p}
\le NC_{w}^N\|w\chi_{B_1}\|_{\p},
\end{equation}
since $P_{j_1}\subset B_1$.

Analogously, we obtain
\begin{equation}\label{dup_winv}
\|w^{-1}\chi_{B_0}\|_{\pri}\lesssim\|w^{-1}\chi_{B_1}\|_{\pri}.
\end{equation}

Finally, if $k$ is the least integer such that $2^k\ge\beta$ and we denote by $\widetilde{B}_0=B(x_0,\frac{r}{2^k})$, it is clear that $\widetilde{B}_0\in\mathcal{B}_\rho$. Then, repeating several times, \eqref{dup_w} and \eqref{dup_winv}, we have
\begin{equation*}
\|w\chi_{B_0}\|_{\p}\|w^{-1}\chi_{B_0}\|_{\pri}\lesssim\|w\chi_{\widetilde{B}_0}\|_{\p}\|w^{-1}\chi_{\widetilde{B}_0}\|_{\pri}\lesssim|\widetilde{B}_0|\lesssim|B_0|.
\end{equation*}

\end{proof}

\begin{rem}
Notice that, even when all the classes $A_{\p}^{\beta\rho,\loc}$ are the same, in general, the membership constant increases with $\beta$.
\end{rem}

\begin{prop}\label{prop:ext}
Let $p\in\mathcal{P}(\R^d)$ and $Q=B(x_0,\be\rho(x_0))$, with $\be>1$. If $w\in A_{\p}^{\rho,\loc}$, then $\widetilde{w}=w\chi_{Q}\in A_{\p}(Q)$. Moreover, the constant $[w]_{A_{\p}(Q)}$ depends only on the constants in \eqref{est rho}, $\be$ and  $[w]_{A_{\p}^{\rho,\loc}}$.
\end{prop}
\begin{proof}
Let $B=B(x,r)$ a ball in $\R^d$ such that $B\cap Q\neq\emptyset$. We need to estimate
\begin{equation*}
\frac1{|B|}\|\chi_{B} \widetilde{w}\|_{\p}\|\chi_{B} \widetilde{w}^{-1}\|_{\pri}.
\end{equation*}
	
	
We distinguish two cases. We start by assuming that $|Q|\le|B|$. In that case,
\begin{equation*}
\begin{split}
\|\chi_{B} \widetilde{w}\|_{\p}\|\chi_{B} \widetilde{w}^{-1}\|_{\pri}&=\|\chi_{B}\chi_{Q}w\|_{\p}\|\chi_{B}\chi_{Q}{w}^{-1}\|_{\pri}\\
&\le\|\chi_Qw\|_{\p}\|\chi_Q{w}^{-1}\|_{\pri}\\
&\le[w]_{ A_{\p}^{\beta\rho,\loc}}|Q|\\
&\le[w]_{ A_{\p}^{\beta\rho,\loc}}|B|,
\end{split}
\end{equation*}
and the constant $[w]_{ A_{\p}^{\beta\rho,\loc}}$ depends only on the constants in \eqref{est rho}, $\be$ and  $[w]_{A_{\p}^{\rho,\loc}}$ as demonstrated in Proposition~\ref{prop: Aploc beta}.
	
Suppose now that $|Q|>|B|$, then $\be\rho(x_0)>r$. Let $y\in B\cap Q$. From \eqref{est rho} to obtain
\begin{equation}\label{ec1}
\rho(x_0)\le c_\rho\rho(y)\left(1+\frac{|y-x_0|}{\rho(x_0)}\right)^{N_0}\le c_{\rho}(1+\be)^{N_0}\rho(y),
\end{equation}
and so
\begin{equation}\label{ec2}
\frac{\rho(x_0)}{\rho(y)}\le c_{\rho}(1+\be)^{N_0}.
\end{equation}

Taking into account~\eqref{ec2} and the fact that $\be\rho(x_0)>r$, from \eqref{est rho} to obtain
\begin{equation}\label{ec3}
\begin{split}
\rho(y)\le c_\rho\rho(x)\left(1+\frac{|y-x|}{\rho(y)}\right)^{N_0}\le c_{\rho}\rho(x)\left(1+\frac{\be\rho(x_0)}{\rho(y)}\right)^{N_0}\\
\le c_{\rho}\big(1+c_{\rho}\be(1+\be)^{N_0}\big)^{N_0}\rho(x).
\end{split}
\end{equation}

Then, by~\eqref{ec1} and~\eqref{ec3},  it follows that $r<\be\rho(x_0)\le\gamma\rho(x)$, with $\gamma=c_{\rho}^2\be(1+\be)^{N_0}\big(1+c_{\rho}\be(1+\be)^{N_0}\big)^{N_0}$, and so $B\in\mathcal{B}_{\gamma\rho}$.

Therefore, from the  Proposition~\ref{prop: Aploc beta}, $w\in A_{\p}^{\gamma\rho,\loc}$ and $[w]_{A_{\p}^{\gamma\rho,\loc}}\le C_{w,\gamma}$, where $C_{w,\gamma}$ only depends on $\gamma$ and $[w]_{A_{\p}^{\rho,\loc}}$, and in this particular case, $C_{w,\gamma}$ only hinges on $\rho$, $\be$ and $[w]_{A_{\p}^{\rho,\loc}}$.

Finally,
\begin{equation*}
\begin{split}
\|\chi_{B}\widetilde{w}\|_{\p}\|\chi_{B} \widetilde{w}^{-1}\|_{\pri}&=\|\chi_{B}\chi_{Q}w\|_{\p}\|\chi_{B}\chi_{Q}w^{-1}\|_{\pri}\\
&\le\|\chi_Bw\|_{\p}\|\chi_B{w}^{-1}\|_{\pri}\\
&\le[w]_{A_{\p}^{\gamma\rho,\loc}}|B|\\
&\lesssim C_{w,\gamma}|B|.
\end{split}
\end{equation*}

Accordingly, $w\in A_{\p}(Q)$ and the constant $[w]_{A_{\p}(Q)}$ depends only of of the constants in \eqref{est rho}, $\be$ and  $[w]_{A_{\p}^{\rho,\loc}}$.

\end{proof}

\section{Main Theorems}

We are now in a position to state one of our main theorems related to weighted norm inequalities for $M^\loc_\rho$.

\begin{thm}\label{teo: acot Mloc}
Let $p\in\mathcal{P}^{\log}(\R^d)$ with $p^{-}>1$. Then, $M_\rho^\loc$ is bounded on $L^{\p}(w)$ if and only if $w\in A_{\p}^{\rho,\loc}$.
\end{thm}
\begin{proof}
Let us note that, following standard arguments, the condition $A_{\p}^{\rho,\loc}$ is nece-ssary for the boundedness of $M_\rho^\loc$ in $L^{\p}(w)$.	
	
Let us now consider $w\in A_{\p}^{\rho,\loc}$. Let $\{Q_j\}_j$ be covering of critical balls given by Proposition~\ref{propo: cubr} and set $\tilde{Q}_j=\be Q_j$, with $\be=c_\rho2^{\frac{N_0}{N_0+1}}+1$.

We observe that for any $x\in\R^d$ and $r\le\rho(x)$, there is $k$ such that $x\in Q_k$, then $B(x,r)\subset B(x,\rho(x))\subset\tilde{Q_k}$ since by~\eqref{est rho}, $\rho(x)\le c_\rho2^{\frac{N_0}{N_0+1}}\rho(x_k)$.

Consequently, for all $f\in L^1_{\loc}(\R^d)$, we have
\begin{equation*}
M_\rho^\loc f(x)\le M( f \chi_{\tilde{Q_k}})(x)
\le\sum_j \chi_{\widetilde{Q}_j}(x)M( f \chi_{\widetilde{Q}_j})(x).
\end{equation*} 

By Proposition~\ref{prop:ext}, we conclude that for each $j$,  $w_j=w\chi_{\widetilde{Q}_j}\in A_{\p}(\widetilde{Q}_j)$ with a constant independent of $j$.

Hence, from Theorem~\ref{teo: local-global} (twice) and Theorem~\ref{thm: Acot M-w} (and Remark~\ref{rem:loc}), we obtain
\begin{equation*}
\begin{split}
\|M_\rho^\loc\|_{\p,w}&\lesssim\bigg\|\sum_{j}\chi_{\tilde{Q}_j}M(f\chi_{\tilde{Q}_j})\bigg\|_{\p,w}\\
&\lesssim\bigg\|\sum_{j}\chi_{\tilde{Q}_j}\frac{\|M(f\chi_{\tilde{Q}_j})\chi_{\tilde{Q}_j}\|_{\p,w}}{\|\chi_{\tilde{Q}_j}\|_{\p,w}}\bigg\|_{\p,w}\\
&=\bigg\|\sum_{j}\chi_{\tilde{Q}_j}\frac{\|M(f\chi_{\tilde{Q}_j})\|_{\p,w_j}}{\|\chi_{\tilde{Q}_j}\|_{\p,w}}\bigg\|_{\p,w}\\
&\lesssim\bigg\|\sum_{j}\chi_{\tilde{Q}_j}\frac{\|f\chi_{\tilde{Q}_j}\|_{\p,w_j}}{\|\chi_{\tilde{Q}_j}\|_{\p,w}}\bigg\|_{\p,w}\\
&=\bigg\|\sum_{j}\chi_{\tilde{Q}_j}\frac{\|f\chi_{\tilde{Q}_j}\|_{\p,w}}{\|\chi_{\tilde{Q}_j}\|_{\p,w}}\bigg\|_{\p,w}\\
&\lesssim\|f\|_{\p,w}.
\end{split}
\end{equation*}
\end{proof}

\begin{rem}
While this paper was in preparation we found the article~\cite{N-S}. In it, the authors defined local weights and determined the weighted inequality for local Hardy--Littlewood maximal operator on the Lebesgue spaces with variable exponent. 

In relation to that, our approach has two differences. On the one hand, in our case the notion of locality is given by the critical radius function $\rho$, while they consider cubes with $|Q|\le1$ (i.e. $\rho\equiv1$). On the other hand, both in the definition of $L^{\p}(w)$ and in that of the  weights $A_{\p}^{\rho,\loc}$, they treat $w\,dx$ as a measure and not as a multiplier.

That is why the methods used there turn out to be quite different from ours.
\end{rem}

\

Our second result yields weighted inequalities concerning the operators  $M_{\rho}^\te$.

\begin{thm}\label{Mtheta}
Let $p\in\mathcal{P}^{\log}(\R^d)$ with $p^{-}>1$. Then,  a weight $w\in A_{\p}^\rho$ if and only if there exists  $\te>0$ such that $M^\te_{\rho}$ is bounded on $L^{\p}(w)$.
\end{thm}
\begin{proof}
First suppose that 	 there exists  $\te>0$ such that $M^\te_{\rho}$ is bounded over $L^{\p}(w)$.

Fix a ball $B=B(x_0,r)$ with $x_0\in\R^d$ and $r>0$. For $x\in B$, using~\eqref{est rho} we have
\begin{equation*}
\rho(x_0)\le c_{\rho}\rho(x)\left(1+\frac{|x-x_0|}{\rho(x_0)}\right)^{N_0}<c_{\rho}\rho(x)\left(1+\frac{2r}{\rho(x_0)}\right)^{N_0}.
\end{equation*}

Then, for any $x\in B$,
\begin{equation}\label{cambio centro}
1+\frac{2r}{\rho(x)}\le1+ c_{\rho}\frac{2r}{\rho(x_0)}\left(1+\frac{2r}{\rho(x_0)}\right)^{N_0}
\le c_\rho\left(1+\frac{2r}{\rho(x_0)}\right)^{N_0+1}.
\end{equation}

Given $\eta>0$, we define the function $A^\eta_Bf$ as
\begin{equation*}
A^\eta_Bf(x)=\frac{\psi_\eta(B)}{|B|}\int_B|f(y)|\,dy\ \chi_B(x)
\end{equation*}
with $\psi_\eta(B)=\big(1+\frac{r}{\rho(x_0)}\big)^{-\eta}$.

Then, taking into account~\eqref{cambio centro}, for $\eta=\te(N_0+1)$ and any $x\in\R^d$, we have
\begin{align*}
A^\eta_Bf(x)&=\frac{\psi_\eta(B)}{|B|}\int_B|f(y)|\,dy\ \chi_B(x)\\
&\le2^{\eta}c_\rho^\te\bigg(1+\frac {2r}{\rho(x)}\bigg)^{-\te}\frac1{|B|}\int_B|f(y)|\,dy\\
&\le2^{d+\eta}c_\rho^\te\bigg(1+\frac {2r}{\rho(x)}\bigg)^{-\te}\frac1{|B(x,2r)|}\int_{B(x,2r)}|f(y)|\,dy\\
&\le2^{d+\eta}c_\rho^\te M_\rho^\te f(x).
\end{align*}

Hence, by Lemma~\ref{lem: formula norma} (exchanging the roles of $\p$ and $\pri$), we obtain
\begin{align*}
\|w\chi_B\|_{\p}\|w^{-1}\chi_B\|_{\pri}&\le2\|w\chi_B\|_{\p}\sup_{\|g\|_{\p}\le1}\int_{B}|g(y)|w^{-1}(y)\,dy\\
&\le2\sup_{\|g\|_{\p}\le1}\left\|w\chi_B\int_{B}|g|w^{-1}\,\right\|_{\p}\\
&=2(\psi_\eta(B))^{-1}|B|\sup_{\|g\|_{\p}\le1}\left\|w\,\chi_B\frac{\psi_\eta(B)}{|B|}\int_{B}|g|w^{-1}\,\right\|_{\p}\\
&=2(\psi_\eta(B))^{-1}|B|\sup_{\|g\|_{\p}\le1}\|A^\eta_B(gw^{-1})w\|_{\p}\\
&\le2^{d+\eta+1}c_\rho^\te\,(\psi_\eta(B))^{-1}|B|\sup_{\|g\|_{\p}\le1}\|M_\rho^\te(gw^{-1})w\|_{\p}\\
&\lesssim(\psi_\eta(B))^{-1}|B|\sup_{\|g\|_{\p}\le1}\|g\|_{\p}\\
&\le\left(1+\frac{r}{\rho(x_0)}\right)^{\eta}|B|,
\end{align*}
which implies that $w\in A_{\p}^\rho$.
	
Suppose now that $w\in A_{\p}^\rho$. Let $\te\ge 0$, which will be specified later. Observe that 
\begin{align*}
M^\te_{\rho}f&\leq M^\te_{\rho,1} f + M^\te_{\rho,2} f,
\end{align*}
where
\begin{equation*}
M^\te_{\rho,1}f(x) = \sup_{r\leq\rho(x)}\bigg(1+\frac r{\rho(x)}\bigg)^{-\te} \frac{1}{|B(x,r)|}\int_{B(x,r)}|f|,
\end{equation*}
and
\begin{equation*}
M^\te_{\rho,2}f(x) = \sup_{r>\rho(x)}\bigg(1+\frac r{\rho(x)}\bigg)^{-\te} \frac{1}{|B(x,r)|} \int_{B(x,r)}|f|.
\end{equation*}

Thereby, we have to check that $M^\te_{\rho,1}$ and $M^\te_{\rho,2}$ are bounded on $L^{\p}(w)$.

Since for every $\te\ge 0$, $M^\te_{\rho,1}\le M_\rho^\loc$ and $A_{\p}^{\rho}\subset A_{\p}^{\rho,\text{loc}}$ the boundedness of $M^\te_{\rho,1}$ follows from Theorem~\ref{teo: acot Mloc}.

Let $\{B_k\}_{k\geq1}$ be a covering provided by  Proposition~\ref{propo: cubr}. Now, for $x\in B_k=B(x_k,\rho(x_k))$ we call $R_j=\{r:2^{j-1}\rho(x)<r\leq2^j\rho(x)\}$, and then we use \eqref{est rho} to obtain
\begin{equation}\label{M2 primera}
\begin{split}
M_{\rho,2}^\te f(x)& \le \sup_{j\geq1} \sup_{r\in R_j}\bigg(1+\frac r{\rho(x)}\bigg)^{-\te}\frac1{|B(x,r)|}\int_{B(x,r)}|f(y)|\,dy\\
&\leq 2^{d+\te} \sup_{j\geq1}\frac{2^{-j(\te+d)}}{|B(x,\rho(x))|}\int_{B(x,2^j\rho(x))}|f(y)|\,dy\\	
& \leq c\sup_{j\geq1} \frac{2^{-j(\te+d)}}{|B_k|}\int_{c_jB_k}|f(y)|\,dy,
\end{split}
\end{equation}
with $c=c_\rho^d 2^{d(1+N_0)+\te}$ and $c_j=2^{j}c_\rho 2^{\frac{N_0}{N_0+1}}+1$.

Moreover, since $w\in A_{\p}^\rho$, there exist $ \sigma$ and $C$
such that
\begin{equation}\label{peso sigma}
\|\chi_B\|_{\p,w}\|\chi_B\|_{\pri,w^{-1}}\le C|B|\left(1+\frac{r}{\rho(x)}\right)^  \sigma,
\end{equation}	
for every ball $B=B(x,r)\subset\R^d$.
	
We denote $B_k^j=c_jB_k$. From H\"older inequality, \eqref{peso sigma} and~\eqref{M2 primera}, it follows
\begin{equation}
\begin{split}
M_{\rho,2}^\te f(x)&\lesssim \sum_{k\ge1}\chi_{B_k}\sup_{j\geq1} \frac{2^{-j(\te+d)}}{|B_k|}\int_{B_k^j}|f(y)|\,dy\\
&\lesssim\sum_{k\ge1}\chi_{B_k}\sup_{j\geq1} \frac{2^{-j(\te+d)}}{|B_k|}\|\chi_{B_k^j}f\|_{\p,w}\|\chi_{B_k^j}\|_{\pri,w^{-1}}\\
&\lesssim\sum_{k\ge1}\chi_{B_k}\sup_{j\geq1} \frac{2^{-j(\te+d)}}{|B_k|}\frac{\|\chi_{B_k^j}f\|_{\p,w}}{\|\chi_{B_k^j}\|_{\p,w}}|B_k^j|c_j^ \sigma\\
&\lesssim\sum_{k\ge1}\chi_{B_k}\sup_{j\geq1} 2^{-j(\te- \sigma)}\frac{\|\chi_{B_k^j}f\|_{\p,w}}{\|\chi_{B_k^j}\|_{\p,w}}\\
&\lesssim\sum_{j\ge1}2^{-j(\te- \sigma)}\sum_{k\ge1}\chi_{B_k^j} \frac{\|\chi_{B_k^j}f\|_{\p,w}}{\|\chi_{B_k^j}\|_{\p,w}}.
\end{split}
\end{equation}	

Thus, by Theorem~\ref{teo: local-global}, we get
\begin{equation*}
\begin{split}
\|M_{\rho,2}^\te f\|_{\p,w}&\le C\left\|\sum_{j\ge1}2^{-j(\te- \sigma)}\sum_{k\ge1}\chi_{B_k^j} \frac{\|\chi_{B_k^j}f\|_{\p,w}}{\|\chi_{B_k^j}\|_{\p,w}}\right\|_{\p,w}\\
&\le C\sum_{j\ge1}2^{-j(\te- \sigma)} \left\|\sum_{k\ge1}\chi_{B_k^j} \frac{\|\chi_{B_k^j}f\|_{\p,w}}{\|\chi_{B_k^j}\|_{\p,w}}\right\|_{\p,w}\\
&\le C\sum_{j\ge1}2^{-j(\te- \sigma)} c_j^{2N_1}\|f\|_{\p,w}\\
&\le C\|f\|_{\p,w}\sum_{j\ge1}2^{-j(\te- \sigma-2N_1)},
\end{split}
\end{equation*}	
where the last series converges taking $\te> \sigma+2N_1$, with $N_1$ is the constant appearing in Proposition~\ref{propo: cubr} (see also Remark~\ref{n-familias}).
\end{proof}

\section{Applications in a Schr\"odinger setting}

In this section we consider a Schr\"{o}dinger operator in $\R^d$ with $d\ge3$,
\begin{equation*}
\mathcal{L}=-\Delta+V,
\end{equation*}
where $V\ge0$, not identically zero, is a function that satisfies for $q>d/2$, the reverse H\"older inequality
\begin{equation}\label{RHq}
\left(\frac1{|B|}\int_B V(y)^q~dy\right)^{1/q}\leq\frac{C}{|B|}\int_B V(y)~dy,
\end{equation}
for every ball $B\subset\R^d$. The set of functions with the last property is usually denoted by $RH_q$.

\begin{rem}
By H\"older inequality we can get that  $RH_q\subset RH_p$, for $q\ge p>1$. One remarkable feature about the $RH_q$ class is that, if $V\in RH_q$ for some $q>1$, then there exists $\epsilon>0$, which depends only on $d$ and the constant $C$ in~\eqref{RHq}, such that  $V\in RH_{q+\epsilon}$. Therefore, it is equivalent to consider $q>d/2$ or $q\ge d/2$ in our hypotheses.
\end{rem}

For a given potential $V\in RH_q$, with $q > d/2$, as in \cite{Shen}, we consider the auxiliary function $\rho_V$ defined for $x\in\R^d$ as
\begin{equation*}
\rho_V(x)=\sup\left\{r>0:\frac1{r^{d-2}}\int_{B(x,r)}V\leq1\right\}.
\end{equation*}
Under the above conditions on $V$ we have $0<\rho(x)<\infty$. In particular, $\rho(x)\sim1$ with $V=C$ and $\rho(x)\sim\frac{1}{1+|x|}$ with $V=|x|^2$ (see~\cite{Martinez}).

Furthermore, according to \cite[Lemma 1.4]{Shen}, if $V\in RH_{q/2}$ the associated function $\rho_V$ verifies \eqref{est rho}.

Thus, if for $f$ in $L^1_{\text{loc}}(\R^d)$, we consider the operators

\begin{equation}
M^\loc_{\rho_V} f(x) = \sup_{r\le\rho_V(x)}\frac{1}{|B(x,r)|} \int_{B(x,r)}|f(y)|\,dy,
\end{equation}

and
\begin{equation}
M^\te_{\rho_V} f(x) = \sup_{r>0}\bigg(1+\frac r{\rho_V(x)}\bigg)^{-\te} \frac{1}{|B(x,r)|} \int_{B(x,r)}|f(y)|\,dy,
\end{equation}
for $\theta>0$, as a consequence of the Theorem~\ref{teo: acot Mloc} and Theorem~\ref{Mtheta}  we get the following theorems. In them, classes $A_{\p}^{\rho_V,\loc}$ and $A_{\p}^{\rho_V}$ refer to those outlined in~\eqref{Ap_local} and~\eqref{Ap_rho} respectively for $\rho=\rho_V$.

\begin{thm}
Let $p\in\mathcal{P}^{\log}(\R^d)$ with $p^{-}>1$ and $V\in RH_q$ with $q > d/2$. Then, $M_{\rho_V}^\loc$ is bounded on $L^{\p}(w)$ if and only if $w\in A_{\p}^{\rho_V,\loc}$.
\end{thm}

\begin{thm}
Let $p\in\mathcal{P}^{\log}(\R^d)$ with $p^{-}>1$ and $V\in RH_q$ with $q > d/2$. Then,  a weight $w\in A_{\p}^{\rho_V}$ if and only if there exists  $\te>0$ such that $M^\te_{\rho_V}$ is bounded on $L^{\p}(w)$.
\end{thm}
\


\end{document}